\documentclass[12pt,leqno]{amsart}
\usepackage[dvips]{graphicx}
\usepackage{amsfonts}
\usepackage{amssymb}
\usepackage{amsmath}
\usepackage{amscd}
\usepackage{graphicx}
\def\DATE{\today}
\newtheorem{theorem}{Theorem}

\newtheorem{corollary}[theorem]{Corollary}

\newtheorem{proposition}[theorem]{Proposition}

\newcommand\C{\mathbb{C}}

\newcommand\g{\mathfrak{g}}

\newcommand\h{\mathfrak{h}}
\newcommand\n{\mathfrak{n}}
\newcommand\K{\mathbb{K}}
\newcommand\Z{\mathbb{Z}}

\newcommand\lr{\left\{ \begin{array}{l}}
\newcommand{\mm}{\mathfrak m }

\title[Nonalgebraic nilpotent Lie algebra]{Examples of nonalgebraic nilpotent Lie algebra }
\author{Elisabeth Remm. Michel Goze}
\date{18 Chevat 5778}
\address{E.R. Universit\'{e} de Haute Alsace, LMIA, 6 rue des Fr\`{e}res Lumi\`{e}re, 68093 Mulhouse, 
M.G. Ramm Algebra Center, F.68800 Rammersmatt}
\email{elisabeth.remm@uha.fr, goze.rac@gmail.com}

\begin{document}

\maketitle

\begin{abstract}
We describe some examples of non abelian nilpotent Lie algebras which are not algebraic. 
\end{abstract}

\medskip

\noindent{\bf Introduction.} All Lie algebras considered in this paper will be finite dimensional and defined over an algebraically closed fixed  $\K$ of characteristic $0$. An algebraic Lie algebra  is a Lie  algebra. The converse  is  in general not true. The problem of determining  if a given Lie algebra  is algebraic or not  is generally difficult. There are some criteria. For example, if a Lie algebra  coincides with its derived subalgebra, then it is algebraic. It is for example the case for the semisimple Lie algebras. But for the abelian Lie algebras, that's not the case. It exists from the dimension  $1$ non algebraic abelian Lie algebras. But we know the structure of the abelian algebraic Lie algebras. It is not the case for nilpotent Lie algebra.
\section{Classical results on algebraic Lie algebras, \cite{H}}

Let $V$ be a finite dimensional $\K$-vector space. Denote by $End(V)$ the $\K$-vector space of endomorphisms of $V$ and by 
$\mathcal{A}ut(V)$ the automorphism group of $V$. Recall that a subgroup of  $\mathcal{A}ut(V)$
 is an algebraic group $G$  if there exist a set $S$ of polynomial functions on $\mathcal{A}ut(V)$ such that $G=\{ f \in \mathcal{A}ut(V), \forall P \in S, P(f)=0\}.$

 With this definition $\mathcal{A}ut(V)$ is an algebraic group.
 
 \noindent {\bf Remark.} A more general definition of algebraic group can be given considering an algebraic group as an affine variety with multiplication and inverse operation which are morphisms of algebraic varieties. In this framework 
   $\mathcal{A}ut(V)$ appears as a principal open subset of the affine space $A^{n^2}.$
  
  The $\K$-vector space $End(V)$  is provided with a Lie algebra structure considering the bracket $[f,g]=f\circ g-g\circ f.$ We denote this Lie algebra by $gl(V)$.

  Let $G$ be an algebraic group of automorphisms of $V$. 
  The Lie algebra of $G$ is the subalgebra of $gl(V)$ wich is composed of the endomorphisms $X$ satisfying $dP(X)=0$ where $P$ is a polynomial in the ideal $I(G)$ generated by $S$. More generally a subalgebra of $gl(V)$ is algebraic if it is the Lie algebra of an algebraic group of automorphisms of $V$. It is difficult to recognize among the Lie algebras those who are algebraic. For example an  abelian Lie algebra is not always algebraic.  In fact,  if we consider a $n$-dimensional algebraic torus, that is a product of $n$ copies of $\K^*$, then its algebraic Lie algebra is diagonalizable, that is it is commutative and all its elements are semi-simple. The algebraic Lie subalgebra of algebraic subgroups of a torus are the linear subspaces of the algebraic Lie algebra of the torus which are defined on the field of rational numbers. 
  
  We recall also some structural properties of an algebraic Lie algebra.  Let $X$ be an endomorphism of $V$. Let $\g(X)$ be the smallest algebraic Lie subalgebra of $gl(V)$ which contains $X$. The elements of $\g(X)$ are called the replica of $X$. Since $X$ generates a one-dimensional Lie subalgebra of $gl(V)$, then $\g(X)$ is the algebraic Lie subalgebra generated by this Lie algebra. We have the following classical result
  \begin{theorem}
A  Lie subalgebra $\g$ of $gl(V)$ is algebraic if and only if for any   $X \in \g$ we have $\g(X) \subset \g$. In particular, if $\g$ is algebraic, the semisimple and nilpotent components $X_s$ and $X_n$ of $X \in \g$ given by the Chevalley-Jordan decomposition of $X$ are also elements of $\g$.
\end{theorem}
 To compute $\g(X)$, we determine in a first step $\g(X_s)$.  Since $X_s$ is semisimple (and $\K=\C)$,  there exists a basis $\{e_1,\cdots,e_n\}$ of $V$ of eigenvectors of $X_s$. Let $\{\lambda_1,\cdots,\lambda_n\}$ the corresponding eigenvalues. The subset $A$ of $gl(V)$ whose elements $Y$ satisfy $Y(e_i)=\mu_ie_i$ is an abelian Lie algebra of $gl(V)$. Now we consider the set $\Lambda$ of $n$-uples of integers  $(p_1,\cdots, p_n)$, $p_I \in \Z$, such that 
$$p_1 \lambda_1 + p_2 \lambda_2 = \cdots +p_n \lambda_n =0.$$
Then, if we denote by $Y_{\mu_1,\cdots,\mu_n}$ the elements of $A$ such that $Y(e_i)=\mu_ie_i$, 
$$\g(X_s)=\{ Y_{\mu_1,\cdots,\mu_n} \in A, \  p_1\mu_1+\cdots +p_n\mu_n=0, \ \ (p_1,\cdots, p_n)\in \Lambda \}.$$
The construction is similar for $\g(X_n)$. The minimal algebraic group $G(X_n)$ containing the unipotent operator $Exp X_n$ consists of all unipotent automorphisms of the form $Exp(tX_n)$. Then $\g(X_n)$ is a one-dimensional abelian algebraic group, isomorphic to the additif abelian group $\K$. 

Recall also the classical decomposition of an algebraic nilpotent Lie algebra.
\begin{proposition}
Let $\g$ be an algebraic nilpotent Lie subalgebra  of $gl(V)$. Let $\frak{n}$ be the set of nilpotent elements of $\g$. Then $\frak{n}$ is an ideal of $\g$ and
\begin{eqnarray}
\label{n+a}
\g= \frak{n} \oplus \frak{a} 
\end{eqnarray}
where $\frak{a}$ is an abelian algebraic Lie subalgebra of $\g$ contained in the center of $\g$ and whose elements are the semisimple elements of $\g$. 
\end{proposition}

\section{An exemple of $3$-dimensional nilpotent Lie algebra which is not algebraic}

We consider the $3$-dimensional linear subspace $\h$ of $gl(4,\C)$ 
$$\left\{ 
\left(
\begin{array}{llll}
x_1+x_2 & x_1+x_2 & 0 & x_1 \\
x_1+x_2 & x_1+x_2 & 0 & x_2 \\
\alpha x_1+(\beta -1) x_2 & \beta x_1+(\alpha +1) x_2 & 0 & x_3 \\
0 & 0 & 0 & 0 
\end{array}
\right), \ x_1,x_2,x_3 \in \C \right\}$$
where  $\alpha, \beta$ are fixed elements of $\C$ such that $\alpha+\beta \neq 0$. A basis of $\h$ is given by the three elements  
$$X_1=
\left(
\begin{array}{llll}
1 & 1 & 0 & 1 \\
1 & 1 & 0 & 0 \\
\alpha  & \beta  & 0 & 0 \\
0 & 0 & 0 & 0 
\end{array}
\right)
,  \ X_2=
\left(
\begin{array}{llll}
1 & 1 & 0 & 0 \\
1 & 1 & 0 & 1 \\
\beta -1  & \alpha +1 & 0 & 0 \\
0 & 0 & 0 & 0 
\end{array}
\right) \ {\rm and } \ \ X_3=\left(
\begin{array}{llll}
0 & 0 & 0 & 0 \\
0 & 0 & 0 & 0 \\
0  & 0 & 0 & 1 \\
0 & 0 & 0 & 0 
\end{array}
\right).
$$
We verify that
$$[X_1,X_2]=X_3, \ [X_1,X_3]=[X_2,X_3]=0.$$
Then $\h$ is a  $3$-dimensional Lie subalgebra of  $gl(4,\C)$. It is isomorphic to the $3$-dimensional Heisenberg Lie algebra, but this isomorphism is not the linear part of isomorphism of algebraic groups. 
Let us consider $X_1$. Its Chevalley-Jordan decomposition  
$$X_1=X_{1,s}+X_{1,n}$$
is given by
 $$X_{1,s}=
\left(
\begin{array}{llll}
1 & 1 & 0 & 1/2 \\
1 & 1 & 0 & 1/2 \\
\frac{\alpha+\beta}{2}  &\frac{\alpha+\beta}{2}   & 0 &\frac{\alpha+\beta}{4} \\
0 & 0 & 0 & 0 
\end{array}
\right), \ X_{1,n}=
\left(
\begin{array}{llll}
0 & 0 & 0 & 1/2 \\
0 & 0 & 0 & -1/2 \\
\frac{\alpha-\beta}{2}   & -\frac{\alpha-\beta}{2}  &0 & -\frac{\alpha+\beta}{4} \\
0 & 0 & 0 & 0 
\end{array}
\right)
$$
Since $X_{1,s} \notin \h$ and also $X_{1,n} \notin \h$, we deduce
\begin{theorem}
The $3$-dimensional nilpotent Lie subalgebra $\h$ of $gl(4,\C)$ is not algebraic.
\end{theorem}

A matrix $X \in \h$ is nilpotent if and only if  $x_1+x_2=0$. In fact the eigenvalues of $X$ are $0$ and $2(x_1+x_2)$. We deduce that the set of nilpotent  matrices of $\h$ is the subspace generated by  $X_1-X_2,X_3$ and it is an abelian ideal of dimension $2$. Let us note that any non trivial matrix of $\h$ is not diagonalisable and  $\h$ doesn't admit a decomposition (\ref{n+a}).

\medskip

We know that any Lie algebra $\g_0$ generates an algebraic Lie algebra $\g_1$ which is the smallest algebraic algebra containing $\g_0$ and these two algebras have the same derived Lie algebra. Let us determinate the algebraic Lie algebra generated by $\h.$ This algebra contains the semi-simple part of any $X \in \h.$ If $X=x_1 X_1+ x_2 X_2+ x_3 X_3,$ 
 $$\displaystyle X_{s}=
\left(
\begin{array}{llll}
x_1+x_2 & x_1+x_2 & 0 & \frac{1}{2}(x_1+x_2)  \\
x_1+x_2 & x_1+x_2 & 0 &(x_1+x_2)  \frac{1}{2} \\
\frac{\alpha+\beta}{2} (x_1+x_2 ) & \frac{\alpha+\beta}{2}(x_1+x_2 )   & 0 &\frac{\alpha+\beta}{4} (x_1+x_2 ) \\
0 & 0 & 0 & 0 
\end{array}
\right)$$
Consider $X_4=\left(
\begin{array}{llll}
1 & 1 & 0 & 1/2 \\
1 & 1 & 0 & 1/2 \\
\frac{\alpha+\beta}{2}  &\frac{\alpha+\beta}{2}   & 0 &\frac{\alpha+\beta}{4} \\
0 & 0 & 0 & 0 
\end{array}
\right)$. Then we have 
$[X_i,X_4]=0$ for $i=1,2,3$
and $\mm=\h \oplus \K\{X_4\}$ is a 4-dimensional Lie algebra, containing $\h$. Moreover, for any $X$ in 
$\mm$, $X_s$ and $X_n$ the semi-simple and nilpotent parts of the Jordan decomposition of $X$ are  in $\mm$.
This Lie algebra is 
$$\mm= \left\{ 
\begin{array}{l}
\left(
\begin{array}{llll}
x_1+x_2+x_4 & x_1+x_2+x_4 & 0 & x_1+\frac{x_4}{2} \\
x_1+x_2+x_4 & x_1+x_2+x_4 & 0 & x_2 +\frac{x_4}{2}\\
\alpha x_1+(\beta -1) x_2+\frac{\alpha+\beta}{2} x_4 & \beta x_1+(\alpha +1) x_2+ \frac{\alpha+\beta}{2} x_4 & 0 & x_3+\frac{\alpha+\beta}{4}x_4  \\
0 & 0 & 0 & 0 
\end{array}
\right) \\
 x_1,x_2,x_3,x_4 \in \C
\end{array}
\right\}, $$
with   fixed elements $\alpha, \beta$ of $\C$ satisfying $\alpha + \beta \neq 0$. Let us note that an element of $\mm$ is nilpotent if and only if $x_1+x_2+x_4=0.$ Then the set $\n_1$ of nilpotent elements of $\mm$ is the $3$-dimensional linear subspace of $\mm$
$$\n_1= \left\{ 
\left(
\begin{array}{llll}
0 & 0 & 0 & x_1/2-x_2/2 \\
0& 0 & 0 & -x_1 /2+x_2/2\\
\frac{\alpha-\beta}{2}  x_1+\frac{-\alpha+\beta-2}{2} x_2 & \frac{-\alpha+\beta}{2} x_1+\frac{\alpha-\beta+2}{2} x_2+ & 0 & x_3+\frac{\alpha+\beta}{4}(-x_1-x_2)  \\
0 & 0 & 0 & 0 
\end{array}
\right)
\right\},$$
The set of  diagonalisable elements is the $1$-dimensional subalgebra
$$\frak{a}_1= \left\{ 
\left(
\begin{array}{llll}
y & y & 0 & y/2 \\
y & y & 0 & y/2\\
\frac{\alpha+\beta}{2} y & \frac{\alpha+\beta}{2} y & 0 & \frac{\alpha+\beta}{4}y  \\
0 & 0 & 0 & 0 
\end{array}
\right)
\right\},$$
with $y \in \C$ and  
$$\mm= \n_1 \oplus \frak{a}_1.$$
Moreover, for any $X \in \mm$, its components $X_s$ and $X_n$ are also in $\mm$.  To study the algebraicity of $\mm$ we have to compute for any $X \in \mm$, the algebraic Lie algebra $\g(X)$ generated by $X$. Let $X_s$ be its semisimple component. The eigenvalues are $0$  which is a  triple root and  the simple root $ 2(x_1+x_2+x_4 )$. We assume that $x_1+x_2+x_4  \neq 0$. The set $\Lambda$ is constituted of $4$-uples of integers $(p_1,p_2,p_3,0)$. Let be $Y$  a semisimple element of $\mm$ commuting with $X_s$ whose  eigenvalues $(\mu_1,\mu_2,\mu_3,\mu_4)$ satisfy $p_1\mu_1+p_2\mu_2+p_3\mu_3=0$ for any $p_1,p_2,p_3 \in \Z$. Then $\mu_1=\mu_2=\mu_3=0$. We deduce that
$$Y=
\left(
\begin{array}{cccc}
m & m & 0 & \frac{m}{2} \\
m &m& 0 & \frac{m}{2} \\
 \frac{m (\alpha +\beta )}{2 } &  \frac{m (\alpha +\beta )}{2 } & 0 &  \frac{m (\alpha +\beta )}{4 } \\
 0 & 0 & 0 & 0 \\
\end{array}
\right)
$$
and $Y$ is an element of $\mm$ corresponding to $x_1=x_2=0,x_4=m.$ Then $\g(X_s) \subset \mm$. Let  $X_n=X-X_s$ be the nilpotent component of $X$. Then 
$$X_n=\left(
\begin{array}{llll}
0 & 0 & 0 & \frac{x_1-x_2}{2} \\
0 & 0 & 0 & \frac{-x_1+x_2}{2} \\
\frac{\alpha-\beta}{2} x_1+\frac{-\alpha+\beta -2}{2}x_2   & \frac{-\alpha+\beta}{2} x_1+\frac{\alpha-\beta+2}{2}x_2  & 0 &x_3-(x_1+x_2 )\frac{\alpha+\beta}{4} \\
0 & 0 & 0 & 0 
\end{array}
\right)
$$
Such vector belongs to the algebraic nilpotent $3$-dimensional Lie algebra $\n_1$ whose elements are  all nilpotent. Recall that a nilpotent subalgebra of $gl(V)$ for some vector space $V$ whose elements are all nilpotent is called unipotent and it is algebraic. But $\g(X_n)$ is contained in this algebra, it is also contained in $\mm$. Then we have

\begin{proposition}
The Lie algebra $\mm$ is algebraic. It is the algebraic Lie algebra generated by the Lie algebra $\h$.
\end{proposition}

\medskip

\noindent{\bf Remark.} The Lie algebra $\h$ is described in \cite{RG1,R1} to proved that there are filiform Lie algebras provided with non complete affine structures.

\section{Other examples in any dimension}
Let $\g$ be a $n$-dimensional nilpotent Lie algebra. Its nilindex is smaller or equal to $n-1$ and $\g$ is called filiform if its nilindex is equal to $n-1$. For example, consider the Lie algebra given in a basis $\{X_1,\cdots,X_n\}$ by 
$$[X_1,X_i]=X_{i+1}, \ i=2,\cdots,n-1, \ \ [X_i,X_j]=0, \ i,j \geq 2.$$
This Lie algebra is commonly called a model of $n$-dimensional filiform Lie algebra and denoted by $L_n$.  

Consider for $n \geq 4$ the $n$-dimensional Lie algebra defined by
$$\scriptsize X_1=\left(
\begin{array}{ccccccccc}
a & a & 0 & 0 & \cdots &0 & 0 & 0 &1\\
a & a & 0 & 0 & \cdots &0 & 0 & 0 &0\\
0 & 0 & 0 & 0 & \cdots &0 & 0 & 0 &0\\
0 & 0 & \frac{1}{2} & 0 & \cdots &0 & 0 & 0 &0\\
0 & 0 & 0 & \frac{2}{3} & \cdots &0 & 0 & 0 &0\\ 
\cdot & \cdot & \cdot & \cdot & \cdots &\frac{n-4}{n-3}& \cdot & \cdot&\cdot\\
\alpha & \beta & 0 & 0 & \cdots &0 & \frac{n-3}{n-2} & 0 &0\\
0 & 0 & 0 & 0 & \cdots &0 & 0 & 0 &0\\
\end{array}
\right) X_2=\left(
\begin{array}{ccccccccc}
a & a & 0 & 0 & \cdots &0 & 0 & 0 &0\\
a & a & 0 & 0 & \cdots &0 & 0 & 0 &1\\
-1 & 1 & 0 & 0 & \cdots &0 & 0 & 0 &0\\
0 & 0 & \frac{1}{2} & 0 & \cdots &0 & 0 & 0 &0\\
0 & 0 & 0 & \frac{1}{3} & \cdots &0 & 0 & 0 &0\\ 
\cdot & \cdot & \cdot & \cdot & \cdots &\frac{n-4}{n-3}& \cdot & \cdot&\cdot\\
\beta & \alpha & 0 & 0 & \cdots &0 & \frac{1}{n-2} & 0 &0\\
0 & 0 & 0 & 0 & \cdots &0 & 0 & 0 &0\\
\end{array}
\right)
$$
with $a \neq 0$ and by
$$X_i=[X_1,X_{i-1}], \ i=2,\cdots,n-1.$$
We deduce that
$$[X_i,X_j]=0, \ i,j \geq 2,$$
and the matrices $X_i$ generate a Lie subalgebra $\h_n$ of $gl(n+1,\C)$ which is isomorphic to $L_n$. The eigenvalues of $X_1$ are $0$ and $2a \neq 0$. Its minimal polynomial is $T^{n-2}(T-2a)$. We deduce that the semisimple component of $X_1$ is 
$$X_{1,s}=\frac{1}{(2a)^{n-3}}X_1^{n-2}.$$
It easily to see that, since $n >3$, this matrix does't belong to $\h_n$. Then $\h_n$ is not an algebraic nilpotent Lie algebra. 
\begin{corollary}
For any $n \geq 3$, there exists a $n$-dimensional non algebraic nilpotent Lie algebra.
\end{corollary}

\medskip

\noindent{\bf Remark.} In some works, one uses a more elusive definition of algebraic Lie algebra considering that a Lie algebra is algebraic if it is isomorphic (with an Lie algebra isomorphism and not a algebraic Lie algebra isomorphism) to an algebraic Lie algebra. From this perspective, any nilpotent lie algebra is algebraic. It is also the case for the abelian Lie algebras, and the study of abelian algebraic groups has not much interest...


\begin{thebibliography}{99}

\bibitem{H} Humphreys, James E. {\it Linear algebraic groups}. Graduate Texts in Mathematics, No. 21. Springer-Verlag, New York-Heidelberg, 1975. xiv+247 pp.


\bibitem{RG1} Remm, E.; Goze, Michel Noncomplete affine structures on Lie algebras of maximal class. Int. J. Math. Math. Sci. 29 (2002), no. 2, 71-77.


\bibitem{R1} Remm, Elisabeth Vinberg algebras associated to some nilpotent Lie algebras. Non-associative algebra and its applications, 347-364, Lect. Notes Pure Appl. Math., 246, Chapman - Hall/CRC, Boca Raton, FL, 2006.
 
\end{thebibliography}
   \end{document}